\documentclass[a4paper,11pt]{amsart}

\usepackage{graphicx}
\usepackage{mathptmx}
\usepackage{amsmath}
\usepackage{amssymb}
\usepackage{enumitem}
\usepackage{xcolor} 
\usepackage{pgfplots}

\newmuskip\pFqmuskip

\newcommand*\pFq[6][8]{%
  \begingroup 
  \pFqmuskip=#1mu\relax
  \mathcode`=\string"8000
  \begingroup\lccode`\~=`\,
  \lowercase{\endgroup\let~}\pFqcomma
  F^{#2}_{#3}{\left(\genfrac..{0pt}{}{#4}{#5}\bigg|#6\right)}%
  \endgroup
}
\newcommand{\pFqcomma}{\mskip\pFqmuskip}

\newtheorem{theorem}{Theorem}[section]

\newtheorem{corollary}[theorem]{Corollary}
\newtheorem{proposition}[theorem]{Proposition}

\begin{document}

\title[]{Degenerate Sheffer-type polynomials and degenerate Sheffer polynomials associated with a random variable}

\author{Taekyun  Kim}
\address{Department of Mathematics, Kwangwoon University, Seoul 139-701, Republic of Korea}
\email{tkkim@kw.ac.kr}
\author{Dae San  Kim}
\address{Department of Mathematics, Sogang University, Seoul 121-742, Republic of Korea}
\email{dskim@sogang.ac.kr}


\begin{abstract}
This paper has two primary contributions. First, we explore degenerate Sheffer-type polynomials, a hybrid of higher-order degenerate Bernoulli and Euler polynomials, and derive their properties. Second, assuming that the moment generating function of $Y$ exists in a neighborhood of the origin, we introduce the degenerate Sheffer polynomials associated with $Y$. We then investigate their properties in general and for the specific cases of uniform and Bernoulli random variables. We also present new results for the higher-order degenerate Bernoulli and Euler polynomials.
\end{abstract}

\maketitle

\markboth{\centerline{\scriptsize Degenerate Sheffer-type polynomials and degenerate Sheffer polynomials}}
{\centerline{\scriptsize Taekyun Kim and Dae San Kim}}

\section{Introduction} 
This paper explores two main areas. First, we delve into degenerate Sheffer-type polynomials, a hybrid of higher-order degenerate Bernoulli and Euler polynomials. Second, we examine degenerate Sheffer polynomials associated with a random variable $Y$, which are a degenerate form of the Appell polynomials associated with $Y$, as previously investigated in \cite{2,20,22}. It's worth noting that the study of various degenerate versions of special polynomials and numbers, originated by Carlitz's work in \cite{4,5}, has seen a renewed surge of interest among mathematicians in recent years (see \cite{3,7,10,13,15,16,23}). We assume that the moment generating function of $Y$ exists in a neighborhood of the origin (see \eqref{23-1}). \par
The paper is structured as follows: \\
$\bullet$ Section 1 reviews degenerate exponentials, degenerate Bernoulli polynomials, and degenerate Euler polynomials. We also briefly remind the reader about uniform and Bernoulli random variables. \\
$\bullet$ Section 2 investigates higher-order degenerate Bernoulli and Euler polynomials, along with Sheffer-type polynomials. After presenting several properties of the higher-order degenerate Bernoulli and Euler polynomials (see Proposition 2.1, Theorems 2.3, 2.4, 2.7), we introduce the degenerate Sheffer-type polynomials $T_{n,\lambda}^{(a,b)}(x)$ (see \eqref{18}). These polynomials represent a hybrid form of the aforementioned higher-order polynomials (see Theorems 2.5, 2.6, 2.8). \\
$\bullet$ Section 3 introduces the degenerate Sheffer polynomials associated with $Y$, denoted as $S_{n,\lambda}^{Y}(x)$ (see \eqref{24}). We derive several properties for these polynomials: \\
\indent $\bullet$ In a general context (see Theorems 3.1, 3.2). \\
\indent$\bullet$ For a uniform random variable on $[0,1]$ (see Theorems 3.3, 3.5, 3.6).\\
\indent$\bullet$ For a Bernoulli random variable with the success probability $\frac{1}{2}$ (see Theorems 3.4, 3.7). \par
Furthermore, Theorem 3.8 establishes an identity for $T_{n,\lambda}^{(a,b)}(x)$ by utilizing the Bernoulli random variable with the success probability of $\frac{1}{2}$. Building on the identity in Theorem 3.8 and other related identities, we deduce an identity for $T_{n,\lambda}^{(a,b)}(x)$ in Theorem 3.9, and for the higher-order degenerate Bernoulli and Euler polynomials in Theorems 3.10 and 3.11. For general references pertaining to this paper, please consult \cite{1,6,14,18,19}.

\vspace{0.1in}

For any nonzero $\lambda\in\mathbb{R}$, the degenerate exponentials are defined by (see \cite{7,10,13,15,16,23})
\begin{equation}
e_{\lambda}^{x}(t)=\sum_{k=0}^{\infty}(x)_{k,\lambda}\frac{t^{k}}{k!},\quad e_{\lambda}(t)=e_{\lambda}^{1}(t), \label{1}
\end{equation}
where 
\begin{equation*}
(x)_{0,\lambda}=1,\quad (x)_{n,\lambda}=x(x-\lambda)(x-2\lambda)\cdots\big(x-(n-1)\lambda\big),\ (n\ge 1). 
\end{equation*}
Note from \eqref{1} that 
\begin{displaymath}
\lim_{\lambda\rightarrow 0}e_{\lambda}^{x}(t)=e^{xt}. 
\end{displaymath}
We also observe that
\begin{equation}
(x+y)_{n,\lambda}=\sum_{k=0}^{n}\binom{n}{k}(x)_{k,\lambda}(y)_{n-k,\lambda}.\label{2}
\end{equation} \par
The degenerate Bernoulli polynomials are defined by (see \cite{13,15})
\begin{equation}
\frac{t}{e_{\lambda}(t)-1}e_{\lambda}^{x}(t)=\sum_{n=0}^{\infty}\beta_{n,\lambda}(x)\frac{t^{n}}{n!}. \label{3}
\end{equation}
When $x=0,\ \beta_{n,\lambda}=\beta_{n,\lambda}(0)$ are called the degenerate Bernoulli numbers. \\
From \eqref{3}, we note that 
\begin{equation}
\beta_{n,\lambda}(x)=\sum_{k=0}^{n}\binom{n}{k}\beta_{k,\lambda}(x)_{n-k,\lambda},\quad (n\ge 0). \label{4}	
\end{equation}
The degenerate Euler polynomials are given by (see \cite{13,15})
\begin{equation}
\frac{2}{e_{\lambda}(t)+1}e_{\lambda}^{x}(t)=\sum_{n=0}^{\infty}\mathcal{E}_{n,\lambda}(x)\frac{t^{n}}{n!}. \label{5}
\end{equation}
When $x=0,\ \mathcal{E}_{n,\lambda}=\mathcal{E}_{n,\lambda}(0),\ (n\ge 0)$, are called the degenerate Euler numbers. \\
By \eqref{5}, we get 
\begin{equation}
\mathcal{E}_{n,\lambda}(x)=\sum_{k=0}^{n}\binom{n}{k}\mathcal{E}_{k,\lambda}(x)_{n-k,\lambda},\quad (n\ge 0). \label{6}
\end{equation}
Note from \eqref{3} and \eqref{5} that 
\begin{displaymath}
\lim_{\lambda\rightarrow 0}\beta_{n,\lambda}(x)=B_{n}(x),\quad \lim_{\lambda\rightarrow 0}\mathcal{E}_{n,\lambda}(x)=E_{n}(x),
\end{displaymath}
where $B_{n}(x)$ and $E_{n}(x)$ are respectively Bernoulli polynomials and Euler polynomials, respectively given by (see \cite{6,12,18})
\begin{equation*}
\frac{t}{e^{t}-1}e^{xt}=\sum_{n=0}^{\infty}B_{n}(x)\frac{t^{n}}{n!},\quad \frac{2}{e^{t}+1}e^{xt}=\sum_{n=0}^{\infty}E_{n}(x)\frac{t^{n}}{n!}. 
\end{equation*} \par
We recall that $Y$ is the uniform random variable on $[a,b]$, denoted by $Y\sim U(a,b)$, if the probability density function of $Y$ is given by (see \cite{14,19})
\begin{equation*}
f_{Y}(x)=\left\{\begin{array}{cc}
\frac{1}{b-a}, & \textrm{if $x\in [a,b]$}, \\
0, & \textrm{otherwise.}
\end{array}\right.
\end{equation*}
Note that 
\begin{displaymath}
E\big[g(Y)\big]=\int_{-\infty}^{\infty}g(y)f_{Y}(y)dy=\frac{1}{b-a}\int_{a}^{b}g(y)dy,
\end{displaymath}
where $g$ is a real-valued function. \\
Also, $Y$ is the Bernoulli random variable with probability of success $p$, denoted by $Y\sim\mathrm{Ber}(p)$, if the probability mass function of $Y$ is given by (see \cite{14,19})
\begin{displaymath}
p(0)=p\{Y=0\}=1-p,\quad p(1)=p\{Y=1\}=p.
\end{displaymath}
Note that 
\begin{displaymath}
E\big[g(Y)\big]=(1-p)g(0)+pg(1). 
\end{displaymath}

\section{Degenerate Sheffer-type polynomials} 
By \eqref{2} and \eqref{4}, we have
\begin{equation}
\beta_{n,\lambda}(x)=\sum_{k=0}^{n}\binom{n}{k}\beta_{k,\lambda}(x)_{n-k,\lambda}=(\beta+x)_{n,\lambda},\quad (n\ge 0), \label{7}	
\end{equation}
with the understanding of replacing $(\beta)_{n,\lambda}$ by $\beta_{n,\lambda}$.
The degenerate Bernoulli numbers satisfy the recurrence relation (see \eqref{3}, \eqref{7}) as follows: 
\begin{equation}
(\beta+1)_{n,\lambda}-\beta_{n,\lambda}=\left\{\begin{array}{cc}
1, & \textrm{if $n=1$} \\
0, & \textrm{otherwise}
\end{array}\right., \quad \beta_{0,\lambda}=1.\label{8}
\end{equation} \\
From the recurrence relation \eqref{8}, we derive
\begin{align*}
&\beta_{0,\lambda}=1,\ \beta_{1,\lambda}=-\frac{1}{2}+\frac{1}{2}\lambda,\ \beta_{2,\lambda}=\frac{1}{6}-\frac{1}{6}\lambda^{2},\ \beta_{3,\lambda}=-\frac{1}{4}\lambda+\frac{1}{4}\lambda^{3}, \\
&\beta_{4,\lambda}=-\frac{1}{30}+\frac{2}{3}\lambda^2-\frac{19}{30}\lambda^{4}, \ \beta_{5,\lambda}=\frac{1}{4}\lambda-\frac{5}{2}\lambda^{3}+\frac{9}{4}\lambda^{5}, \dots.
\end{align*} \par
From \eqref{2} and \eqref{6}, we note that 
\begin{equation}
\mathcal{E}_{n,\lambda}(x)=\sum_{k=0}^{n}\binom{n}{k}\mathcal{E}_{k,\lambda}(x)_{n-k,\lambda}=\big(\mathcal{E}+x\big)_{n,\lambda},\quad (n\ge 0), \label{9}	
\end{equation}
with the understanding of replacing $(\mathcal{E})_{n,\lambda}$ by $\mathcal{E}_{n,\lambda}$.
The degenerate Euler numbers satisfy the recurrence relation (see \eqref{5}, \eqref{9}) as follows: 
\begin{equation}
\big(\mathcal{E}+1\big)_{n,\lambda}+\mathcal{E}_{n,\lambda}=\left\{\begin{array}{cc}
2, & \textrm{if $n=0$} \\
0, & \textrm{if $n\ne 0$}
\end{array}\right.,\quad \mathcal{E}_{0,\lambda}=1. \label{10}
\end{equation} \\
By using the recurrence relation \eqref{10}, we obtain
\begin{align*}
&\mathcal{E}_{0,\lambda}=1,\ \mathcal{E}_{1,\lambda}=-\frac{1}{2},\ \mathcal{E}_{2,\lambda}=\frac{1}{2}\lambda,\ \mathcal{E}_{3,\lambda}=\frac{1}{4}-\lambda^{2}, \\
&\mathcal{E}_{4,\lambda}= -\frac{3}{2}\lambda + 3 \lambda^{3}, \ \mathcal{E}_{5,\lambda}=-\frac{1}{2}+\frac{35}{4}\lambda^{2}-12\lambda^{4}, \dots.
\end{align*}

For $\alpha\in\mathbb{R}$, the higher-order degenerate Bernoulli polynomials are defined by 
\begin{equation}
\bigg(\frac{t}{e_{\lambda}(t)-1}\bigg)^{\alpha}e_{\lambda}^{x}(t)=\sum_{n=0}^{\infty}\beta_{n,\lambda}^{(\alpha)}(x)\frac{t^{n}}{n!}. \label{11}
\end{equation}
When $x=0,\ \beta_{n,\lambda}^{(\alpha)}=\beta_{n,\lambda}^{(\alpha)}(0)$ are called the higher-order degenerate Bernoulli numbers. \\
The higher-order degenerate Euler polynomials are given by 
\begin{equation}
\bigg(\frac{2}{e_{\lambda}(t)+1}\bigg)^{\alpha}e_{\lambda}^{x}(t)=\sum_{n=0}^{\infty}\mathcal{E}_{n,\lambda}^{(\alpha)}(x)\frac{t^{n}}{n!}. \label{12}
\end{equation}
When $x=0,\ \mathcal{E}_{n,\lambda}^{(\alpha)}=\mathcal{E}_{n,\lambda}^{(\alpha)}(0)$ are called the higher-order degenerate Euler numbers. \\
Note that 
\begin{displaymath}
\lim_{\lambda\rightarrow 0}B_{n,\lambda}^{(\alpha)}(x)=B_{n}^{(\alpha)}(x),\quad \lim_{\lambda\rightarrow 0}\mathcal{E}_{n,\lambda}^{(\alpha)}(x)=E_{n}^{(\alpha)}(x), 
\end{displaymath}
where $B_{n}^{(\alpha)}(x)$ and $E_{n}^{(\alpha)}(x)$ are the higher-order Bernoulli and higher-order Euler polynomial, respectively given by (see \cite{18})
\begin{equation*}
\bigg(\frac{t}{e^{t}-1}\bigg)^{\alpha}e^{xt}=\sum_{n=0}^{\infty}B_{n}^{(\alpha)}(x)\frac{t^{n}}{n!},\quad \bigg(\frac{2}{e^{t}+1}\bigg)^{\alpha}e^{xt}=\sum_{n=0}^{\infty}E_{n}^{(\alpha)}(x)\frac{t^{n}}{n!}. 
\end{equation*}
From \eqref{11} and \eqref{12}, we note that 
\begin{equation}
\beta_{n,\lambda}^{(0)}(x)=(x)_{n,\lambda}\quad \mathrm{and}\quad \mathcal{E}_{n,\lambda}^{(0)}(x)=(x)_{n,\lambda},\quad (n\ge 0). \label{13}	
\end{equation}

The following proposition is immediate from the definitions \eqref{11} and \eqref{12}.
\begin{proposition}
For $a,b\in\mathbb{R}$, we have 
\begin{equation*}
\beta_{n,\lambda}^{(a+b)}(x+y)=\sum_{k=0}^{n}\binom{n}{k}\beta_{k,\lambda}^{(a)}(x)\beta_{n-k,\lambda}^{(b)}(y),
\end{equation*}
and 
\begin{equation*}
\mathcal{E}_{n,\lambda}^{(a+b)}(x+y)=\sum_{k=0}^{n}\binom{n}{k}\mathcal{E}_{k,\lambda}^{(a)}(x)\mathcal{E}_{n-k,\lambda}^{(b)}(y). 
\end{equation*}
\end{proposition}
From Proposition 2.1 and \eqref{13}, we obtain the following corollary.
\begin{corollary}
For $a\in\mathbb{R}$, we have 
\begin{displaymath}
\beta_{n,\lambda}^{(a)}(x+y)=\sum_{k=0}^{n}\binom{n}{k}\beta_{k,\lambda}^{(a)}(x)(y)_{n-k,\lambda},
\end{displaymath}
and 
\begin{displaymath}
\mathcal{E}_{n,\lambda}^{(a)}(x+y)=\sum_{k=0}^{n}\binom{n}{k}\mathcal{E}_{k,\lambda}^{(a)}(x)(y)_{n-k,\lambda}. 
\end{displaymath}
\end{corollary}
Now, we observe that 
\begin{align}
&\sum_{n=0}^{\infty}\big(\beta_{n,\lambda}^{(\alpha)}(x+1)-\beta_{n,\lambda}^{(\alpha)}(x)\big)\frac{t^{n}}{n!} =t\bigg(\frac{t}{e_{\lambda}(t)-1}\bigg)^{\alpha-1}e_{\lambda}^{x}(t) \label{14}\\
&=\sum_{n=0}^{\infty}(n+1)\beta_{n,\lambda}^{(\alpha-1)}(x)\frac{t^{n+1}}{(n+1)!}=\sum_{n=0}^{\infty}n\beta_{n-1,\lambda}^{(\alpha-1)}(x)\frac{t^{n}}{n!},\nonumber
\end{align}
and 
\begin{align}
\sum_{n=0}^{\infty}\Big(\mathcal{E}_{n,\lambda}^{(\alpha)}(x+1)+\mathcal{E}_{n,\lambda}^{(\alpha)}(x)\Big)\frac{t^{n}}{n!}&=2\bigg(\frac{2}{e_{\lambda}(t)+1}\bigg)^{\alpha-1}e_{\lambda}^{x}(t)\label{15}\\
&=\sum_{n=0}^{\infty}2\mathcal{E}_{n,\lambda}^{(\alpha-1)}(x)\frac{t^{n}}{n!}. \nonumber
\end{align}
Therefore, by \eqref{14} and \eqref{15}, we obtain the following theorem. 
\begin{theorem}
For $n\ge 0$, and  $\alpha\in\mathbb{R}$, we have 
\begin{displaymath}
\beta_{n,\lambda}^{(\alpha)}(x+1)-\beta_{n,\lambda}^{(\alpha)}(x)=n\beta_{n-1,\lambda}^{(\alpha-1)}(x),
\end{displaymath}
and 
\begin{displaymath}
\mathcal{E}_{n,\lambda}^{(\alpha)}(x+1)+\mathcal{E}_{n,\lambda}^{(\alpha)}(x)=2\mathcal{E}_{n,\lambda}^{(\alpha-1)}(x). 
\end{displaymath}
\end{theorem}
From \eqref{3} and \eqref{5}, we note that 
\begin{align}
&\sum_{n=0}^{\infty}\beta_{n,\lambda}(x)\frac{t^{n}}{n!}=\frac{t}{e_{\lambda}(t)-1}e_{\lambda}^{x}(t)=\frac{1}{2}\frac{t}{e_{\lambda}(t)-1}\frac{2}{e_{\lambda}(t)+1}\big(e_{\lambda}(t)+1\big)e_{\lambda}^{x}(t) \label{16} \\
&=\frac{1}{2}\frac{t}{e_{\lambda}(t)-1}e_{\lambda}(t)\frac{2}{e_{\lambda}(t)+1}e_{\lambda}^{x}(t)+\frac{1}{2}\frac{t}{e_{\lambda}(t)-1}\frac{2}{e_{\lambda}(t)+1}e_{\lambda}^{x}(t)\nonumber\\
&=\frac{1}{2}\sum_{n=0}^{\infty}\sum_{k=0}^{n}\binom{n}{k}\Big(\beta_{k,\lambda}(1)+\beta_{k,\lambda}\Big)\mathcal{E}_{n-k,\lambda}(x)\frac{t^{n}}{n!}. \nonumber
\end{align}
Observe that $\beta_{k,\lambda}(1)=\beta_{k,\lambda}+\delta_{k,1}$ (see \eqref{8}), where $\delta_{k,1}$ is the Kronecker's delta. By \eqref{16}, we get 
\begin{align}
\beta_{n,\lambda}(x)&=\frac{1}{2}\sum_{k=0}^{n}\binom{n}{k}\Big(2\beta_{k,\lambda}+\delta_{k,1}\Big)\mathcal{E}_{n-k,\lambda}(x) \label{17}	\\
&=\frac{n}{2}\mathcal{E}_{n-1,\lambda}(x)+\sum_{k=0}^{n}\binom{n}{k}\beta_{k,\lambda}\mathcal{E}_{n-k,\lambda}(x). \nonumber
\end{align}
Therefore, by \eqref{17}, we obtain the following theorem. 
\begin{theorem}
For $n\ge 0$, we have 
\begin{displaymath}
\beta_{n,\lambda}(x)= \frac{n}{2}\mathcal{E}_{n-1,\lambda}(x)+\sum_{k=0}^{n}\binom{n}{k}\beta_{k,\lambda}\mathcal{E}_{n-k,\lambda}(x).
\end{displaymath}
\end{theorem}
Now, we consider the {\it{degenerate Sheffer-type polynomials}} given by   
\begin{equation}
\bigg(\frac{t}{e_{\lambda}(t)-1}\bigg)^{a}\bigg(\frac{2}{e_{\lambda}(t)+1}\bigg)^{b}e_{\lambda}^{x}(t)=\sum_{n=0}^{\infty}T_{n,\lambda}^{(a,b)}(x)\frac{t^{n}}{n!},\label{18}
\end{equation}
where $a,b\in\mathbb{R}$. \par 
The following proposition follows immediately from the definition \eqref{18}.
\begin{proposition}
For $n\ge 0$, and $a,b\in\mathbb{R}$, we have 
\begin{displaymath}
T_{n,\lambda}^{(a,b)}(x+y)=\sum_{k=0}^{n}\binom{n}{k}\beta_{k,\lambda}^{(a)}(x)\mathcal{E}_{n-k,\lambda}^{(b)}(y),
\end{displaymath}
and 
\begin{displaymath}
T_{n,\lambda}^{(a,b)}(x)=\sum_{k=0}^{n}\binom{n}{k}T_{k,\lambda}^{(a,b)}(0)(x)_{n-k,\lambda}. 
\end{displaymath}
\end{proposition}
From \eqref{18} and \eqref{3}, we have 
\begin{align}
\sum_{n=0}^{\infty}T_{n,\lambda}^{(a,b)}(x)\frac{t^{n}}{n!}&=\bigg(\frac{t}{e_{\lambda}(t)-1}\bigg)^{a-1}\bigg(\frac{2}{e_{\lambda}(t)+1}\bigg)^{b}\bigg(\frac{t}{e_{\lambda}(t)-1}\bigg)e_{\lambda}^{x}(t) \label{19}\\
&=\sum_{k=0}^{\infty}T_{k,\lambda}^{(a-1,b)}(0)\frac{t^{k}}{k!} \sum_{l=0}^{\infty}\beta_{l,\lambda}(x)\frac{t^{l}}{l!}\nonumber\\
&=\sum_{n=0}^{\infty}\sum_{k=0}^{n}\binom{n}{k}T_{k,\lambda}^{(a-1,b)}(0)\beta_{n-k,\lambda}(x)\frac{t^{n}}{n!}.\nonumber	
\end{align}
Similarly to \eqref{19}, by \eqref{18} and \eqref{5}, we get 
\begin{align}
\sum_{n=0}^{\infty}T_{n,\lambda}^{(a,b)}(x)\frac{t^{n}}{n!}&=\bigg(\frac{t}{e_{\lambda}(t)-1}\bigg)^{a}\bigg(\frac{2}{e_{\lambda}(t)+1}\bigg)^{b-1} \frac{2}{e_{\lambda}(t)+1} e_{\lambda}^{x}(t)\label{20} \\
&=\sum_{k=0}^{\infty}T_{k,\lambda}^{(a,b-1)}(0)\frac{t^{k}}{k!}\sum_{l=0}^{\infty}\mathcal{E}_{l,\lambda}(x)\frac{t^{l}}{l!}\nonumber\\
&=\sum_{n=0}^{\infty}\sum_{k=0}^{n}\binom{n}{k}T_{k,\lambda}^{(a,b-1)}(0)\mathcal{E}_{n-k,\lambda}(x)\frac{t^{n}}{n!}. \nonumber
\end{align}
Therefore, by \eqref{19} and \eqref{20}, we obtain the following theorem. 
\begin{theorem}
For $n \ge 0$, and $a,b\in\mathbb{R}$, we have 
\begin{align*}
T_{n,\lambda}^{(a,b)}(x)&=\sum_{k=0}^{n}\binom{n}{k}T_{k,\lambda}^{(a-1,b)}(0)\beta_{n-k,\lambda}(x)\\
&=\sum_{k=0}^{n}\binom{n}{k}T_{k,\lambda}^{(a,b-1)}(0)\mathcal{E}_{n-k,\lambda}(x). 
\end{align*}
\end{theorem}
Now, we observe that 
\begin{align}
&\sum_{n=0}^{\infty}\bigg(\sum_{k=0}^{n}\binom{n}{k}\beta_{k,\lambda}\mathcal{E}_{n-k,\lambda}(x)\bigg)\frac{t^{n}}{n!}=\bigg(\frac{t}{e_{\lambda}(t)-1}\bigg)\bigg(\frac{2}{e_{\lambda}(t)+1}\bigg)e_{\lambda}^{x}(t)
\label{21}\\
&=\frac{2t}{e_{\lambda}^{2}(t)-1}e_{\lambda}^{x}(t)=\frac{2t}{e_{\frac{\lambda}{2}}(2t)-1}e_{\frac{\lambda}{2}}^{\frac{x}{2}}(2t)=\sum_{n=0}^{\infty}\beta_{n,\frac{\lambda}{2}}\bigg(\frac{x}{2}\bigg)2^{n}\frac{t^{n}}{n!}.\nonumber
\end{align}
Therefore, by comparing the coefficients on both sides of \eqref{21},we obtain the following theorem. 
\begin{theorem}
For $n\ge 0$, we have 
\begin{displaymath}
2^{n}\beta_{n,\frac{\lambda}{2}}\bigg(\frac{x}{2}\bigg)=\sum_{k=0}^{n}\binom{n}{k}\beta_{k,\lambda}\mathcal{E}_{n-k,\lambda}(x). 
\end{displaymath}
\end{theorem}
From \eqref{18}, we note that 
\begin{align}
&\bigg(\frac{t}{e_{\lambda}(t)-1}\bigg)^{a} \bigg(\frac{t}{e_{\lambda}(t)-1}\bigg)^{b}e_{\lambda}^{x+1}(t)-\bigg(\frac{t}{e_{\lambda}(t)-1}\bigg)^{a} \bigg(\frac{t}{e_{\lambda}(t)-1}\bigg)^{b}e_{\lambda}^{x}(t)\label{22}\\
&=t \bigg(\frac{t}{e_{\lambda}(t)-1}\bigg)^{a-1} \bigg(\frac{t}{e_{\lambda}(t)-1}\bigg)^{b}e_{\lambda}^{x}(t)=\sum_{n=0}^{\infty}(n+1)T_{n,\lambda}^{(a-1,b)}(x)\frac{t^{n+1}}{(n+1)!}\nonumber\\
&=\sum_{n=0}^{\infty}nT_{n-1,\lambda}^{(a-1,b)}(x)\frac{t^{n}}{n!}. \nonumber	
\end{align}
Therefore, by \eqref{18} and \eqref{22}, we obtain the following theorem. 
\begin{theorem}
For $n\ge 0$, we have 
\begin{equation}
T_{n,\lambda}^{(a,b)}(x+1)-T_{n,\lambda}^{(a,b)}(x)=nT_{n-1,\lambda}^{(a-1,b)}(x). \label{23}
\end{equation}
\end{theorem}

\section{Degenerate Sheffer polynomials associated with a random variable}
Assume that $Y$ is a random variable whose moment generating function exists in some neighborhood of the origin (see \cite{2,11,12,15,16,17,21,23}), that is,
\begin{equation}
E\Big[e^{Yt}\Big]<\infty,\quad (|t|<r), \quad\textrm{for some $r>0$}. \label{23-1}
\end{equation}
Let $(Y_{j})_{j\ge 1}$ be a sequence of mutually independent copies of the random variable $Y$, and let 
\begin{displaymath}
Y^{(k)}=Y_{1}+\cdots+Y_{k},\ (k\ge 1),\quad\mathrm{with}\quad Y^{(0)}=0.
\end{displaymath} \par
Now, we define the {\it{degenerate Sheffer polynomials associated with $Y$},} $S_{n,\lambda}^{Y}(x),\ (n \ge 0)$, by 
\begin{equation}
\frac{1}{E\big[e_{\lambda}^{Y}(t)\big]}e_{\lambda}^{x}(t)=\sum_{n=0}^{\infty}S_{n,\lambda}^{Y}(x)\frac{t^{n}}{n!}.\label{24} 	
\end{equation}
Note that 
\begin{align}
\sum_{n=0}^{\infty}E\Big[S_{n,\lambda}^{Y}(x+Y)\Big]\frac{t^{n}}{n!}&=E\bigg[\sum_{n=0}^{\infty}S_{n,\lambda}^{Y}(x+Y)\frac{t^{n}}{n!}\bigg]\label{25}	\\
&=E\bigg[\frac{1}{E\big[e_{\lambda}^{Y}(t)\big]}e_{\lambda}^{Y+x}(t)\bigg]= \frac{e_{\lambda}^{x}(t)}{E\big[e_{\lambda}^{Y}(t)\big]}E\big[e_{\lambda}^{Y}(t)\big]\nonumber\\
&=e_{\lambda}^{x}(t)=\sum_{n=0}^{\infty}(x)_{n,\lambda}\frac{t^{n}}{n!}.\nonumber 
\end{align}
Therefore, by comparing the coefficients on both sides of \eqref{25}, we obtain the following theorem. 
\begin{theorem}
For $n\ge 0$, we have 
\begin{displaymath}
E\big[S_{n,\lambda}^{Y}(x+Y)\big]=(x)_{n,\lambda}.
\end{displaymath}
\end{theorem}
Let $Y_{1}$ and $Y_{2}$ be independent random variables. Then we have 
\begin{align}
\sum_{n=0}^{\infty}S_{n,\lambda}^{Y_{1}+Y_{2}}(x+y)\frac{t^{n}}{n!}&=\frac{1}{E\big[e_{\lambda}^{Y_{1}+Y_{2}}(t)\big]}e^{(x+y)t} \label{26}	\\
&=\frac{1}{E\big[e_{\lambda}^{Y_{1}}(t)\big] E\big[e_{\lambda}^{Y_{2}}(t)\big] }e^{xt}e^{yt}\nonumber\\
&=\sum_{k=0}^{\infty}S_{k,\lambda}^{Y_{1}}(x)\frac{t^{k}}{k!}\sum_{j=0}^{\infty}S_{j,\lambda}^{Y_{2}}(y)\frac{t^{j}}{j!} \nonumber\\
&=\sum_{n=0}^{\infty}\sum_{k=0}^{n}\binom{n}{k}S_{k,\lambda}^{Y_{1}}(x)S_{n-k,\lambda}^{Y_{2}}(y)\frac{t^{n}}{n!}. \nonumber
\end{align}
In particular, for $y=0$ and $Y_{2}=0$, we get 
\begin{align}
\sum_{n=0}^{\infty}S_{n,\lambda}^{Y_{1}}(x)\frac{t^{n}}{n!}&=\frac{1}{E\big[e_{\lambda}^{Y_{1}}(t)\big]}e_{\lambda}^{x}(t)\label{27}\\
&=\sum_{k=0}^{\infty}S_{k,\lambda}^{Y_{1}}(0)\frac{t^{k}}{k!}\sum_{l=0}^{\infty}(x)_{l,\lambda}\frac{t^{l}}{l!}\nonumber\\
&=\sum_{n=0}^{\infty}\sum_{k=0}^{n}\binom{n}{k}S_{k,\lambda}^{Y_{1}}(0)(x)_{n-k,\lambda}\frac{t^{n}}{n!}.\nonumber	
\end{align}
Therefore, by \eqref{26} and \eqref{27}, we obtain the following theorem. 
\begin{theorem}
For $n\ge 0$, let $Y_{1},Y_{2}$ be independent random variables. Then we have 
\begin{displaymath}
S_{n,\lambda}^{Y_{1}+Y_{2}}(x+y)=\sum_{k=0}^{n}\binom{n}{k}S_{k,\lambda}^{Y_{1}}(x)S_{n-k,\lambda}^{Y_{2}}(y).
\end{displaymath}
In particular, for $y=0$ and $Y_{2}=0$, we have 
\begin{displaymath}
S_{1,\lambda}^{Y_{1}}(x)=\sum_{k=0}^{n}\binom{n}{k}S_{k,\lambda}^{Y_{1}}(0)(x)_{n-k,\lambda}. 
\end{displaymath}
\end{theorem}
For $Y\sim U[0,1]$,  we have 
\begin{equation}
E\big[e_{\lambda}^{Y}(t)\big]=\int_{0}^{1}e_{\lambda}^{y}(t)dy=\frac{\lambda t}{\log(1+\lambda t)}\frac{1}{t}\big(e_{\lambda}(t)-1\big). \label{28}
\end{equation}
Thus, by \eqref{28}, we have 
\begin{align}
\sum_{n=0}^{\infty}S_{n,\lambda}^{Y}(x)\frac{t^{n}}{n!}&=\frac{1}{E\big[e_{\lambda}^{Y}(t)\big]}e_{\lambda}^{x}(t)=\frac{\log(1+\lambda t)}{\lambda t}\frac{t}{e_{\lambda}(t)-1}e_{\lambda}^{x}(t) \label{29}\\
&=\sum_{l=0}^{\infty}\frac{(-\lambda)^{l}l!}{l+1}\frac{t^{l}}{l!} \sum_{k=0}^{\infty}\beta_{k,\lambda}(x)\frac{t^{k}}{k!}\nonumber\\
&=\sum_{n=0}^{\infty}\sum_{k=0}^{n}\binom{n}{k}\beta_{k,\lambda}(x)\frac{(-\lambda)^{n-k}(n-k)!}{n-k+1}\frac{t^{n}}{n!}.\nonumber
\end{align}
Therefore, by \eqref{29}, we obtain the following theorem, 
\begin{theorem}
For $Y\sim U[0,1]$, we have 
\begin{displaymath}
S_{n,\lambda}^{Y}(x)=\sum_{k=0}^{n}\binom{n}{k}\beta_{k,\lambda}(x)\frac{(-\lambda)^{n-k}(n-k)!}{n-k+1}, \quad (n \ge 0).
\end{displaymath}
\end{theorem}
Let $Y\sim\mathrm{Ber}(1/2)$. Then we have 
\begin{equation}
E\big[e_{\lambda}^{Y}(t)\big]=\frac{1}{2}+\frac{1}{2}e_{\lambda}(t)=\frac{e_{\lambda}(t)+1}{2}.\label{30}
\end{equation}
By \eqref{30}, we get 
\begin{equation}
\sum_{n=0}^{\infty}S_{n,\lambda}^{Y}(x)\frac{t^{n}}{n!}=\frac{1}{E\big[e_{\lambda}^{Y}(t)\big]}e_{\lambda}^{x}(t)=\frac{2}{e_{\lambda}(t)+1}e_{\lambda}^{x}(t)=\sum_{n=0}^{\infty}\mathcal{E}_{n,\lambda}(x)\frac{t^{n}}{n!}. \label{31}
\end{equation}
Therefore, by \eqref{31}, we obtain the following theorem. 
\begin{theorem}
For $Y\sim\mathrm{Ber}(1/2)$, and $n\ge 0$, we have 
\begin{displaymath}
S_{n,\lambda}^{Y}(x)=\mathcal{E}_{n,\lambda}(x). 
\end{displaymath}
\end{theorem}
For $Y\sim U[0,1]$, and $m\in\mathbb{N}$, by \eqref{28}, we have 
\begin{align}
E\Big[e_{\lambda}^{Y^{(m)}}(t)\Big]&=E\Big[e_{\lambda}^{Y_{1}+Y_{2}+\cdots+Y_{m}}(t)\Big]=\big(e_{\lambda}(t)-1\big)^{m}\bigg(\frac{\lambda}{\log(1+\lambda t)}\bigg)^{m} \label{32}\\
&=\bigg(\frac{e_{\lambda}(t)-1}{t}\bigg)^{m}\bigg(\frac{\lambda t}{\log(1+\lambda t)}\bigg)^{m}.\nonumber
\end{align}
By \eqref{32}, we get 
\begin{align}
\sum_{n=0}^{\infty}S_{n,\lambda}^{Y^{(m)}}(x)\frac{t^{n}}{n!}&=\frac{e_{\lambda}^{x}(t)}{E\big[e_{\lambda}^{Y^{(m)}}(t)\big]}=\bigg(\frac{\log(1+\lambda t)}{\lambda t}\bigg)^{m}\bigg(\frac{t}{e_{\lambda}(t)-1}\bigg)^{m}e_{\lambda}^{x}(t)\label{33}\\
&=\frac{m!}{\lambda^{m}}\sum_{j=m}^{\infty}S_{1}(j,m)\lambda^{j}\frac{t^{j-m}}{j!}\sum_{k=0}^{\infty}\beta_{k,\lambda}^{(m)}(x)\frac{t^{k}}{k!} \nonumber \\
&=\sum_{j=0}^{\infty}\frac{S_{1}(j+m,m)\lambda^{j}}{\binom{j+m}{m}}\frac{t^{j}}{j!}\sum_{k=0}^{\infty}\beta_{k,\lambda}^{(m)}(x)\frac{t^{k}}{k!}\nonumber\\
&=\sum_{n=0}^{\infty}\sum_{k=0}^{n}S_{1}(n-k+m,m)\lambda^{n-k}\frac{\binom{n}{k}}{\binom{n-k+m}{m}}\beta_{k,\lambda}^{(m)}(x)\frac{t^{n}}{n!},\nonumber
\end{align}
where $S_{1}(j,m)$ are the Stirling number of the first kind defined by 
\begin{equation*}
\frac{1}{k!}\log^{k}(1+t)=\sum_{n=k}^{\infty}S_{1}(n,k)\frac{t^{n}}{n!},\quad (k\ge 0). 
\end{equation*}
Therefore, by \eqref{33}, we obtain the following theorem. 
\begin{theorem}
For $Y\sim U[0,1]$, and $m,n\ge 0$, we have 
\begin{displaymath}
S_{n,\lambda}^{Y^{(m)}}(x)=\sum_{k=0}^{n}\cfrac{\binom{n}{k}}{\binom{n-k+m}{m}}\lambda^{n-k}S_{1}(n-k+m,m)\beta_{k,\lambda}^{(m)}(x). 
\end{displaymath}
\end{theorem}
Let $m,l$ be integers with $m \ge l \ge 0$. Now, by \eqref{24} and \eqref{25}, we observe that 
\begin{align}
\sum_{n=0}^{\infty}E\Big[S_{n,\lambda}^{Y^{(m)}}\big(x+Y^{(l)}\big)\Big]\frac{t^{n}}{n!}&=\frac{e_{\lambda}^{x}(t)}{E\big[e_{\lambda}^{Y^{(m)}}(t)\big]}E\Big[e_{\lambda}^{Y^{(l)}}(t)\Big]\label{34}	\\
&=\frac{e_{\lambda}^{x}(t)}{E\big[e_{\lambda}^{Y^{(m)}-Y^{(l)}}(t)\big]}=\sum_{n=0}^{\infty}S_{n,\lambda}^{Y^{(m)}-Y^{(l)}}(x)\frac{t^{n}}{n!}. \nonumber 
\end{align}
Thus, by \eqref{34}, we get 
\begin{equation}
E\Big[S_{n,\lambda}^{Y^{(m)}}\big(x+Y^{(l)}\big)\Big]=S_{n,\lambda}^{Y^{(m)}-Y^{(l)}}(x),\quad (n\ge 0, \ m \ge l \ge 0). \label{35}
\end{equation}
Therefore, by \eqref{35} and Theorem 3.5, we obtain the following theorem. 
\begin{theorem}
For $n\ge 0$, and $Y\sim U[0,1]$, we have 
\begin{align*}
&\sum_{k=0}^{n}\frac{\binom{n}{k}}{\binom{n-k+m}{m}}\lambda^{n-k}S_{1}(n-k+m,m)E\bigg[\beta_{n,\lambda}^{(m)}\bigg(x+Y^{(l)}\bigg)\bigg] \\
&=\sum_{k=0}^{n}\frac{\binom{n}{k}}{\binom{n-k+m-l}{m-l}}\lambda^{n-k}S_{1}(n-k+m-l,m-l)\beta_{k,\lambda}^{(m-l)}(x), 
\end{align*}
where $m,l$ are integers with $m \ge l \ge 0$. 
\end{theorem}
For $Y\sim\mathrm{Ber}(1/2)$, we have 
\begin{align}
\sum_{n=0}^{\infty}S_{n,\lambda}^{Y^{(m)}}(x)\frac{t^{n}}{n!}&= \frac{e_{\lambda}^{x}(t)}{E\big[e_{\lambda}^{Y^{(m)}}(t)\big]}=\frac{1}{E\big[e_{\lambda}^{Y_{1}+\cdots+Y_{m}}(t)\big]}e_{\lambda}^{x}(t) \label{36}\\
&=\frac{1}{E\big[e_{\lambda}^{Y_{1}}(t)\big]\cdots E\big[e_{\lambda}^{Y_{m}}(t)\big]}e_{\lambda}^{x}(t)=\bigg(\frac{2}{e_{\lambda}(t)+1}\bigg)^{m}e_{\lambda}^{x}(t) \nonumber\\
&=\sum_{n=0}^{\infty}\mathcal{E}_{n,\lambda}^{(m)}(x)\frac{t^{n}}{n!},\quad (m \ge 0). \nonumber
\end{align}
By \eqref{36}, we get 
\begin{equation}
S_{n,\lambda}^{Y^{(m)}}(x)=\mathcal{E}_{n,\lambda}^{(m)}(x).\label{37}
\end{equation}
Therefore, by \eqref{35} and \eqref{37}, we obtain the following theorem. 
\begin{theorem}
For $Y\sim\mathrm{Ber}(1/2)$, and $n\ge 0$, we have 
\begin{displaymath}
E\bigg[S_{n,\lambda}^{Y^{(m)}}\bigg(x+ Y^{(l)}\bigg)\bigg]=\mathcal{E}_{n,\lambda}^{(m-l)}(x), 
\end{displaymath}
where $m,l$ are integers with $m \ge l \ge 0$. 
\end{theorem}
For $Y\sim\mathrm{Ber}(1/2)$, we have 
\begin{align}
\sum_{n=0}^{\infty}E\Big[T_{n,\lambda}^{(a,b)}(x+Y)\Big]\frac{t^{n}}{n!}&=E\bigg[\sum_{n=0}^{\infty}T_{n,\lambda}^{(a,b)}(x+Y)\frac{t^{n}}{n!}\bigg] \label{38}\\
&=\bigg(\frac{t}{e_{\lambda}(t)-1}\bigg)^{a}\bigg(\frac{2}{e_{\lambda}(t)+1}\bigg)^{b}e_{\lambda}^{x}(t)E\big[e_{\lambda}^{Y}(t)\big] \nonumber	\\ 
&=\bigg(\frac{t}{e_{\lambda}(t)-1}\bigg)^{a}\bigg(\frac{2}{e_{\lambda}(t)+1}\bigg)^{b}e_{\lambda}^{x}(t)\frac{e_{\lambda}(t)+1}{2} \nonumber\\
&=\bigg(\frac{t}{e_{\lambda}(t)-1}\bigg)^{a}\bigg(\frac{2}{e_{\lambda}(t)+1}\bigg)^{b-1}e_{\lambda}^{x}(t)\nonumber\\
&=\sum_{n=0}^{\infty}T_{n,\lambda}^{(a,b-1)}(x)\frac{t^{n}}{n!}. \nonumber
\end{align}
On the other hand, by Theorem 2.5, we get 
\begin{align}
E\big[T_{n,\lambda}^{(a,b)}(x+Y)\big]&=\sum_{k=0}^{n}\binom{n}{k}T_{k,\lambda}^{(m,l)}(0)E\big[(x+Y)_{n-k,\lambda}\big] \label{39} 	\\
&=\frac{1}{2}\sum_{k=0}^{n}\binom{n}{k}T_{k,\lambda}^{(a,b)}(0)\Big((x+1)_{n-k,\lambda}+(x)_{n-k,\lambda}\Big)\nonumber\\
&=\frac{1}{2}\Big(T_{n,\lambda}^{(a,b)}(x+1)+T_{n,\lambda}^{(a,b)}(x)\Big). \nonumber
\end{align}
From \eqref{38} and \eqref{39}, we obtain the following theorem. 
\begin{theorem}
For $n\ge 0$, we have 
\begin{equation}
T_{n,\lambda}^{(a,b-1)}(x)= \frac{1}{2}\Big(T_{n,\lambda}^{(a,b)}(x+1)+T_{n,\lambda}^{(a,b)}(x)\Big),\label{40}
\end{equation}
where $a, b \in \mathbb{R}$.
\end{theorem}
From \eqref{23} and \eqref{40}, we derive the following theorem. 
\begin{theorem}
For $n\ge 0$, we have 
\begin{equation}
T_{n,\lambda}^{(a,b)}(x)=T_{n,\lambda}^{(a,b-1)}(x)-\frac{n}{2}T_{n-1,\lambda}^{(a-1,b)}(x),\label{41}
\end{equation}
where $a,b \in \mathbb{R}$. 
\end{theorem}
By Theorem 2.5 and \eqref{41}, we get 
\begin{align}
&\sum_{k=0}^{n}\binom{n}{k}\beta_{k,\lambda}^{(a)}(x)\mathcal{E}_{n-k,\lambda}^{(b)}(y)=T_{n,\lambda}^{(a,b)}(x+y)=T_{n,\lambda}^{(a,b-1)}(x+y)-\frac{n}{2}T_{n-1}^{(a-1,b)}(x+y) \label{42}\\
&=\sum_{k=0}^{n}\binom{n}{k}\beta_{k,\lambda}^{(a)}(x)\mathcal{E}_{n-k,\lambda}^{(b-1)}(y)-\frac{n}{2}\sum_{k=0}^{n-1}\binom{n-1}{k}\beta_{k,\lambda}^{(a-1)}(x)\mathcal{E}_{n-k-1,\lambda}^{(b)}(y) \nonumber\\ 
&=\sum_{k=0}^{n}\binom{n}{k}\beta_{k,\lambda}^{(a)}(x)\mathcal{E}_{n-k,\lambda}^{(b-1)}(y)-\sum_{k=0}^{n}\binom{n}{k}\frac{k}{2}\beta_{k-1,\lambda}^{(a-1)}(x)\mathcal{E}_{n-k,\lambda}^{(b)}(y).\nonumber
\end{align}
From \eqref{42}, we obtain the following theorem. 
\begin{theorem}
For $n\ge 0$, we have
\begin{equation*}
\sum_{k=0}^{n}\binom{n}{k}\beta_{n,\lambda}^{(a)}(x)\mathcal{E}_{n-k,\lambda}^{(b-1)}(y)=\sum_{k=0}^{n}\binom{n}{k}\bigg[\beta_{n,\lambda}^{(a)}(x)+\frac{k}{2}\beta_{k-1,\lambda}^{(a-1)}(x)\bigg]\mathcal{E}_{n-k,\lambda}^{(b)}(y),
\end{equation*}
where $a, b \in \mathbb{R}$.
\end{theorem}
Let $b=1$ in Theorem 3.10. Then we have 
\begin{equation}
\sum_{k=0}^{n}\binom{n}{k}\beta_{k,\lambda}^{(a)}(x)(y)_{n-k,\lambda}=\sum_{k=0}^{n}\binom{n}{k}\bigg[\beta_{n,\lambda}^{(a)}(x)+\frac{k}{2}\beta_{k-1,\lambda}^{(a-1)}(x)\bigg]\mathcal{E}_{n-k,\lambda}(y). \label{43}
\end{equation}
Thus, by \eqref{43} and Corollary 2.2, we get 
\begin{equation}
\beta_{n,\lambda}^{(a)}(x+y)=\sum_{k=0}^{n}\binom{n}{k}\bigg[\beta_{k,\lambda}^{(a)}(x)+\frac{k}{2}\beta_{k-1,\lambda}^{(a-1)}(x)\bigg]\mathcal{E}_{n-k,\lambda}(y).\label{44}
\end{equation} \par
Now, we observe that 
\begin{equation}
\sum_{k=0}^{n}\binom{n}{k}\frac{k}{2}\beta_{k-1,\lambda}^{(a-1)}(x)\mathcal{E}_{n-k,\lambda}^{(b)}(y)=\frac{n}{2}\sum_{k=0}^{n-1}\binom{n-1}{k}\beta_{k,\lambda}^{(a-1)}(x)\mathcal{E}_{n-k-1,\lambda}^{(b)}(y).\label{45}
\end{equation}
Let $a=1$ in \eqref{45}. Then, by Corollary 2.2, we have 
\begin{align}
\sum_{k=0}^{n}\binom{n}{k}\frac{k}{2}\beta_{k-1,\lambda}^{(0)}(x)\mathcal{E}_{n-k,\lambda}^{(b)}(y)&=\frac{n}{2}\sum_{k=0}^{n-1}\binom{n-1}{k}(x)_{k,\lambda}\mathcal{E}_{n-1-k,\lambda}^{(b)}(y) \label{46} \\
&=\frac{n}{2}\mathcal{E}_{n-1,\lambda}^{(b)}(x+y). \nonumber
\end{align}
By Theorem 3.10 with $a=1$ and \eqref{46}, we get 
\begin{align}
\sum_{k=0}^{n}\binom{n}{k}\beta_{k,\lambda}(x)\mathcal{E}_{n-k,\lambda}^{(b)}(y)&=\sum_{k=0}^{n}\binom{n}{k}\beta_{k,\lambda}(x)\mathcal{E}_{n-k,\lambda}^{(b-1)}(y)-\sum_{k=0}^{n}\binom{n}{k}\frac{k}{2}\beta_{k-1,\lambda}^{(0)}(x)\mathcal{E}_{n-k,\lambda}^{(b)}(y)\label{47}\\
&=\sum_{k=0}^{n}\binom{n}{k}\beta_{k,\lambda}(x)\mathcal{E}_{n-k,\lambda}^{(b-1)}(y)-\frac{n}{2}\mathcal{E}_{n-1,\lambda}^{(b)}(x+y). \nonumber
\end{align}
From \eqref{47}, we note that 
\begin{align}
\frac{n}{2}\mathcal{E}_{n-1,\lambda}^{(b)}(x+y)&=\sum_{k=0}^{n}\binom{n}{k}\beta_{k,\lambda}(x)\Big(\mathcal{E}_{n-k,\lambda}^{(b-1)}(y)-\mathcal{E}_{n-k,\lambda}^{(b)}(y)\Big) \label{48}\\
&=\sum_{k=0}^{n}\binom{n}{k}\beta_{n-k,\lambda}(x)\Big(\mathcal{E}_{k,\lambda}^{(b-1)}(y)-\mathcal{E}_{k,\lambda}^{(b)}(y)\Big),\nonumber	
\end{align}
where $n \ge 0$. \\
Thus, for $n \ge 1$, from \eqref{48} we have
\begin{align}
\mathcal{E}_{n-1,\lambda}^{(b)}(x+y)&=\frac{2}{n}\sum_{k=0}^{n}\binom{n}{k}\beta_{n-k,\lambda}(x)\Big(\mathcal{E}_{k,\lambda}^{(b-1)}(y)-\mathcal{E}_{k,\lambda}^{(b)}(y)\Big) \label{49} \\
&=\sum_{k=1}^{n}\binom{n-1}{k-1}\frac{2}{k}\beta_{n-k,\lambda}(x)\Big(\mathcal{E}_{k,\lambda}^{(b-1)}(y)-\mathcal{E}_{k,\lambda}^{(b)}(y)\Big) \nonumber\\
&=\sum_{k=0}^{n-1}\frac{2}{k+1}\binom{n-1}{k}\beta_{n-1-k,\lambda}(x)\Big(\mathcal{E}_{k+1,\lambda}^{(b-1)}(y)-\mathcal{E}_{k+1,\lambda}^{(b)}(y)\Big). \nonumber
\end{align}
Therefore, by \eqref{44} and \eqref{49}, we obtain the following theorem. 
\begin{theorem}
For $a, b \in \mathbb{R}$ and $n\ge 0$, we have
\begin{align*}
&\beta_{n,\lambda}^{(a)}(x+y)=\sum_{k=0}^{n}\binom{n}{k}\bigg[\beta_{k,\lambda}^{(a)}(x)+\frac{k}{2}\beta_{k-1,\lambda}^{(a-1)}(x)\bigg]\mathcal{E}_{n-k,\lambda}(y), \\
&\mathcal{E}_{n,\lambda}^{(b)}(x+y)=\sum_{k=0}^{n}\binom{n}{k}\frac{2}{k+1}\beta_{n-k,\lambda}(x)\Big(\mathcal{E}_{k+1,\lambda}^{(b-1)}(y)-\mathcal{E}_{k+1,\lambda}^{(b)}(y)\Big). 
\end{align*}
\end{theorem}

\section{Conclusion} 
In this paper, we studied a hybrid of higher-order degenerate Bernoulli and Euler polynomials, namely the degenerate Sheffer-type polynomials. Then, we examined a degenerate form of the Appell polynomials associated with a random variable $Y$, namely the degenerate Sheffer polynomials associated with $Y$. \par
We notice from \eqref{24} and \eqref{30} that, for $Y \sim \mathrm{Ber}(1/2)$, we have
\begin{equation*}
\frac{2}{e_{\lambda}(t)+1}e_{\lambda}^{x}(t)=\sum_{n=0}^{\infty}S_{n,\lambda}^{\mathrm{Ber}(1/2)}(x)\frac{t^{n}}{n!},
\end{equation*}
so that $S_{n,\lambda}^{\mathrm{Ber}(1/2)}(x)=\mathcal{E}_{n,\lambda}(x)$, the degenerate Euler polynomials in \eqref{5}. \\
While, for $Y \sim U[0,1]$, we observe from \eqref{24} and \eqref{28} that 
\begin{equation}
\frac{\log(1+\lambda t)}{\lambda(e_{\lambda}(t)-1)}e_{\lambda}^{x}(t)=\sum_{n=0}^{\infty}S_{n,\lambda}^{U[0,1]}(x)\frac{t^{n}}{n!}, \label{50}
\end{equation}
so that $S_{n,\lambda}^{U[0,1]}(x)$ are the `degenerate Bernoulli polynomials arising from Volkenborn integral' (see \cite{9}). These polynomials are different from the Carlitz degenerate Bernoulli polynomials $\beta_{n,\lambda}(x)$ in \eqref{3} and \eqref{4}. However, as $\lambda \rightarrow 0$, both $\beta_{n,\lambda}(x)$ and $S_{n,\lambda}^{U[0,1]}(x)$ tend to the Bernoulli polynomial $B_{n}(x)$. For details on this, the interested reader may refer to \cite{9}. \par
It is one of our future research projects to continue to explore various aspects and applications of degenerate versions of many special polynomials and numbers.

\end{document}